\documentclass[fleqn]{ifacconf}% Modiefied by Ebihara

\usepackage{graphicx}      % include this line if your document contains figures
\usepackage{natbib}        % required for bibliography

\usepackage{mathtools}
\usepackage{amsfonts}
\usepackage{amssymb}
\usepackage{color}
\usepackage{comment}

\newcommand{\nc}{\newcommand}

\nc{\mb}[1]{\mathbb{#1}}
\nc{\ml}[1]{\mathcal{#1}}
\nc{\mf}[1]{\mathbf{#1}}
\nc{\mr}[1]{\mathrm{#1}}
\nc{\slope}[3]{{#1}\in\mr{slope}[#2,#3]}
\nc{\trace}{\mr{trace}}
\nc{\He}{\mr{He}}
\nc{\rank}[1]{\mr{rank}(#1)}

\nc{\w}{w^{\ast}}
\nc{\z}{z^{\ast}}
\nc{\1}{\mf{1}_m}
\nc{\Md}{M_{\mr{d}}}
\nc{\oMd}{\overline{M}_{\mr{d}}}
\nc{\Mod}{M_{\mr{od}}}
\nc{\oMod}{\overline{M}_{\mr{od}}}

\nc{\slopePhi}[2]{{\mf{\Phi}}_{#1,#2}^{\lowercase{m}}}
\nc{\oddPhi}{\mf{\Phi}_{\mr{\lowercase{odd}}}^m}
\nc{\phiwc}{\phi_{\mr{wc}}}
\nc{\Phiwc}{\Phi_{\mr{wc}}}

\nc{\Pd}[1]{\ml{P}_{\mr{d}}\left( #1 \right)}
\nc{\Pod}[1]{\ml{P}_{\mr{od}}\left( #1 \right)}

\nc{\primalfirst}
  {
  \begin{bmatrix}
    -P+A^TPA & A^TPB \\
    B^TPA & B^TPB      
  \end{bmatrix}
  }
   
\nc{\primalsecond}[3]
  {
  (\ast)^T\ml{M}(#1,#2,#3)\begin{bmatrix}
    C & D \\
    0 & I_m
    \end{bmatrix}
  }    
  
\nc{\matrixY}
  {
  \begin{bmatrix}
    -\mu C & -\mu D + I_m 
  \end{bmatrix} H \begin{bmatrix}
      \nu C^T \\
      \nu D^T - I_m
    \end{bmatrix}
  }  
  
\nc{\RematrixY}
  {
  \begin{bmatrix}
    0_{m,n} & I_m 
  \end{bmatrix} H \begin{bmatrix}
      C^T \\
      D^T - I_m
    \end{bmatrix}
  }    
  
\nc{\dualfirst}
  {
  -H_{11}+\begin{bmatrix}
    A & B
    \end{bmatrix} H \begin{bmatrix}
      A & B
      \end{bmatrix}^T
  }        
\begin{document}
\begin{frontmatter}

\title{%
Detecting Destabilizing Nonlinearities in \\
Absolute Stability Analysis of \\
Discrete-Time Feedback Systems\thanksref{footnoteinfo}} 

\thanks[footnoteinfo]
{
This work was supported by JSPS KAKENHI Grant Number JP21H01354 and
Japan Science and Technology Agency (JST) as part of
Adopting Sustainable Partnerships for Innovative Research Ecosystem (ASPIRE), Grant Number JPMJAP2402. 
This work was also supported by the AI Interdisciplinary Institute ANITI funding, 
through the French "Cluster AI" program under 
the Grant agreement n°ANR-23-IACL-0002 as well as 
by the National Research Foundation, Prime Minister's Oﬃce, 
Singapore under its Campus for Research Excellence and 
Technological Enterprise (CREATE) programme.
}

\author[Kyushu-S]{Hibiki Gyotoku} 
\author[Kyushu-F]{Tsuyoshi Yuno} 
\author[Kyushu-F]{Yoshio Ebihara}
\author[LAAS]{Dimitri Peaucelle}
\author[LAAS]{Sophie Tarbouriech}
\author[LAAS]{Victor Magron}

\address[Kyushu-S]{%
Graduate School of Information Science and Electrical Engineering, Kyushu University, 
Fukuoka 8190395, Japan
(gyotoku.hibiki.327@s.kyushu-u.ac.jp).}
\address[Kyushu-F]{%
Faculty of Information Science and Electrical Engineering, \\Kyushu University, 
Fukuoka 8190395, Japan. }
\address[LAAS]{%
\mbox{LAAS-CNRS, Universit\'{e} de Toulouse, CNRS, F-31400, Toulouse, France.}}

\begin{abstract}                % Abstract of 50--100 words
This paper is concerned with the absolute stability analysis of 
discrete-time feedback systems with slope-restricted nonlinearities. 
By employing static O'Shea-Zames-Falb multipliers in the framework of 
integral quadratic constraints, we can obtain a certificate for the absolute stability 
in the form of a linear matrix inequality (LMI).  
However, since this LMI certificate is only a sufficient condition, 
we cannot draw any definite conclusion if the LMI turns out to be infeasible.  
To address this issue, we focus on the dual LMI that is feasible if and only 
if the original (primal) LMI is infeasible.  
As the main result, if the dual solution satisfies a certain rank condition, 
we prove that we can detect a destabilizing nonlinearity within 
the assumed class of slope-restricted nonlinearities as well as 
a non-zero equilibrium point of the resulting feedback system, 
thereby we can conclude that the system of interest is never absolutely stable. 
The effectiveness of the technical results is demonstrated through numerical examples.  
\end{abstract}

\end{frontmatter}
%===============================================================================

%%%%%%%%%%%%%%%%%%%%%%%%%%%%%%%%%%%%%%%%%%%%%%%%%%%%%%%%%%%%%%%%%%%%%%%%%
%%%%%%%%%%%%%%%%%%%%%%%%%%%%%%%%%%%%%%%%%%%%%%%%%%%%%%%%%%%%%%%%%%%%%%%%%
\section{Introduction}
\label{sec:intro}
In recent years,
 there has been a growing attention to control theoretic methods dealing with 
the analysis of optimization algorithms (\cite{Lessard_SIAM2016}), 
 stability analysis of dynamical neural networks (NNs) 
(\cite{Revay_IEEE2021,Ebihara_ECC2021}), 
 and performance analysis of control systems driven by NNs 
 (\cite{Yin_IEEE2022,Scherer_IEEE2022,Souza_2023}). 
These algorithms, NNs and control systems can be modeled as 
feedback systems with 
linear time-invariant (LTI) systems and nonlinear operators, 
and this allows us to apply control theory for their analysis.  
In particular, since 
NNs have a large number of various nonlinear operators such as 
hyperbolic tangent ($\tanh$), sigmoid, and rectified linear unit (ReLU), 
it is required to establish a solid method for the 
analysis of feedback systems with a variety of nonlinearities. 

Against this background, 
we deal with the absolute stability analysis of 
discrete-time nonlinear feedback systems. 
Here, a feedback system is said to be absolutely stable if 
its origin in the state space is globally asymptotically stable for 
all nonlinear operators belonging to an assumed class (\cite{Khalil_2002}).    
In this regard, this paper focuses on the class of 
slope-restricted and repeated nonlinear operators 
(\cite{Valmorbida_IEEE2018,Fagundes_2024}). 
For this class of nonlinearities, it is well known that  
O'Shea-Zames-Falb (OZF) multipliers 
(\cite{O'Shea_IEEE1967,Zames_SIAM1968,Carrasco_EJC2016})
are effective  in the framework of Integral Quadratic Constraints (IQCs) 
(\cite{Megretski_IEEE1997}). 
This approach enables us to obtain 
linear matrix inequality (LMI) conditions to 
ensure the absolute stability of the feedback systems. 
However, since these LMIs are only sufficient conditions, 
we cannot conclude anything on the absolute stability 
if the LMIs are numerically infeasible.   

To address this issue, 
this paper focuses on the dual LMIs of 
those IQC-based (primal) LMIs, 
where the dual LMI is feasible if and only if the (primal) LMI is infeasible (\cite{Scherer_EJC2006}). 
The main results of this paper can be summarized in the following way; 
if the solution of the dual LMI satisfies a specific rank condition, then 
1) 
it is possible to detect a nonlinear operator
that destabilizes the target system  
within the assumed class of the slope-restricted and repeated nonlinearities,   
2) 
we can also identify a non-zero equilibrium point that proves the instability of 
the target feedback system with the detected destabilizing nonlinear operator,   
3)
and hence we can conclude that the target system is never absolutely stable. 
The validity of the technical results is demonstrated through numerical examples. 

This study is a continuation of 
our recent studies on nonlinear feedback system analysis using dual LMIs.  
For \emph{continuous-time} feedback systems with repeated ReLU nonlinearities, 
\cite{Yuno_IFAC2024} examined the dual LMI of an IQC-based LMI and 
derived a rank condition on the dual solution to conclude that 
the target system is not globally asymptotically stable. 
Building on this, in our recent research (\cite{Gyotoku_ECC2025}), 
we extended the results in (\cite{Yuno_IFAC2024}) to 
the absolute stability analysis of continuous-time feedback systems with 
slope-restricted and repeated nonlinearities.   
This paper is an extension of \cite{Gyotoku_ECC2025} to discrete-time feedback systems.  
The novel aspects of dealing with discrete-time systems
are summarized in Remark~\ref{rem:novel} of Section \ref{sec:slope-odd}.  
We finally note that the contents of subsections~\ref{ssec:OZF}, \ref{ssec:Construction_DHD}, and~\ref{ssec:Construction_DD} can be found in the paper \citep{Gyotoku_ECC2025}, but we reiterate them to make the current paper self-contained.

Notation: \ 
The set of $n{\times}m$ real matrices (resp. with nonnegative entries)
is denoted by $\mb{R}^{n{\times}m}$ (resp. $\mb{R}_+^{n{\times}m}$). 
For a matrix $A\in\mb{R}^{n{\times}m}$, $A_{i,j}$ stands for the $(i,j)$ entry of $A$.
We denote the $n{\times}n$ identity matrix,
 the $n{\times}n$ zero matrix, and the $n{\times}m$ zero matrix by 
 $I_n, \ 0_n, \ 0_{n,m}$, respectively. 
For a matrix $A$, we write $A \ge 0$ to denote that
 $A$ is entrywise nonnegative. 
We denote the set of
 $n{\times}n$ real symmetric, positive semidefinite, and positive definite matrices by 
 $\mb{S}^n, \ \mb{S}_+^n, \ \mb{S}_{++}^n$, respectively. 
For $A\in\mb{S}^n$, we write $A \succ 0$ (resp. $A \prec 0$) to denote that
 $A$ is positive (resp. negative) definite. 
For $A\in\mb{R}^{n{\times}n}$ and $B\in\mb{R}^{n{\times}m}$,
 $(\ast)^TAB$ is a shorthand notation of $B^TAB$. 
For $v\in\mb{R}^n$, we denote by $\|v\|$ its standard Euclidean norm. 
The induced norm of 
 a (possibly nonlinear) operator $\Phi : \mb{R}^m \rightarrow \mb{R}^n$ is defined by 
 $\|\Phi\| := \underset{v\in\mb{R}^m\backslash\{0\}}{\mr{sup}}\dfrac{\|\Phi(v)\|}{\|v\|}$. 
For $A\in\mb{R}^{n{\times}n}$, 
 we define $|A|_{\mr{d}}\in\mb{R}^{n{\times}n}$ by 
 $|A|_{\mr{d},i,i} = A_{i,i}$ and $|A|_{\mr{d},i,j} = -|A_{i,j}| \ (i \ne j)$. 
We also define
 \begin{equation*}
   \begin{aligned}
     &
     \mb{D}^m := \{
       M\in\mb{R}^{m{\times}m} : M_{i,j} = 0 \ (i \ne j, \ i,j = 1,\ldots,m)
       \}, \\
     &
     \mb{OD}^m := \{
       M\in\mb{R}^{m{\times}m} : M_{i,i} = 0 \ (i = 1,\ldots,m)
       \}.  
   \end{aligned}
 \end{equation*}
We finally define the (canonical) projections onto $\mb{D}^m$ and $\mb{OD}^m$ as $\ml{P}_{\mr{d}} : \mb{R}^{m{\times}m} \ni M \mapsto \mathrm{diag}(M_{1,1},M_{2,2},\ldots,M_{m,m})\in\mb{D}^m$ and $\ml{P}_{\mr{od}} : \mb{R}^{m{\times}m} \ni M \mapsto M - \ml{P}_{\mr{d}}(M) \in \mb{OD}^m$, respectively.
%

%
%%%%%%%%%%%%%%%%%%%%%%%%%%%%%%%%%%%%%%%%%%%%%%%%%%%%%%%%%%%%%%%%%%%%%%%%%
%%%%%%%%%%%%%%%%%%%%%%%%%%%%%%%%%%%%%%%%%%%%%%%%%%%%%%%%%%%%%%%%%%%%%%%%%
\section{Basics on IQC and Static O'Shea-Zames-Falb Multipliers}
%%%%%%%%%%%%%%%%%%%%%%%%%%%%%%%%%%%%%%%%
\subsection{Basic Results on IQC-based Stability Analysis}

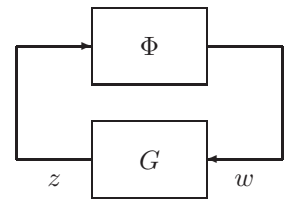
\begin{figure}[b] % block diagram
  \centering
  \begin{picture}(3.5,2.5)(0,0)
\put(0,2){\vector(1,0){1}}
\put(1,1.5){\framebox(1.5,1){$\Phi$}}
\put(2.5,2){\line(1,0){1}}
\put(3.5,2){\line(0,-1){1.5}}
\put(3.5,0.5){\vector(-1,0){1}}
\put(3,0.3){\makebox(0,0)[t]{$w$}}
\put(1,0){\framebox(1.5,1){$G$}}
\put(1,0.5){\line(-1,0){1}}
\put(0.5,0.3){\makebox(0,0)[t]{$z$}}
\put(0,0.5){\line(0,1){1.5}}
\end{picture}
  \caption{Discrete-Time Nonlinear Feedback System $\Sigma$.}
  \label{fig:system} 
\end{figure}
Let us consider the discrete-time feedback system $\Sigma$ shown in 
Fig.~\ref{fig:system}. 
Here, $G$ is an LTI system described by 
 \begin{equation}
   G : \begin{cases}
     x(k+1) &= Ax(k) + Bw(k), \\
     z(k) &= Cx(k) + Dw(k) 
     \end{cases}
 \label{eq:LTI}    
 \end{equation} 
where $A\in\mb{R}^{n{\times}n}$, $B\in\mb{R}^{n{\times}m}$, 
$C\in\mb{R}^{m{\times}n}$, and $D\in\mb{R}^{m{\times}m}$. 
On the other hand, 
$\Phi : \mb{R}^m \rightarrow \mb{R}^m$ stands for a static nonlinear operator described by 
 \begin{equation}
   w(k) = \Phi(z(k)).
 \label{eq:nonlinearity}   
 \end{equation} 
We assume that $A$ is Schur stable, i.e., its spectral radius is less than one,  and  
the feedback system $\Sigma$ is well-posed.
In this paper, we investigate the absolute stability analysis of 
 the feedback system $\Sigma$ in IQC framework with multipliers. 
Here, the feedback system $\Sigma$ is said to be 
stable if  its origin in the state space is globally asymptotically stable, 
and further $\Sigma$ is said to be  absolutely stable if 
it is stable for all nonlinear operators in an assumed class. 
The basic result on IQC-based stability analysis is 
presented in the next proposition. 

%% IQC basic result
\begin{prop} \label{prop:IQC} (\cite{Megretski_IEEE1997}) \\
Let us define the set of multipliers $\mf{\Pi}^{\star} \subset \mb{S}^{2m}$ by 
 \begin{equation*}
   \mf{\Pi}^{\star} := \left\{
     \Pi\in\mb{S}^{2m} : (\ast)^T\Pi\begin{bmatrix}
       \zeta \\
       \Phi(\zeta)
       \end{bmatrix} \ge 0 \quad \forall{\zeta}\in\mb{R}^m
   \right\}.        
 \end{equation*}
Then, the feedback system $\Sigma$ is stable if 
 there exist $P \succ 0$ and $\Pi\in\mf{\Pi}^{\star}$ such that 
 \begin{equation*}
   \primalfirst + (\ast)^T\Pi\begin{bmatrix}
     C & D \\
     0 & I_m
     \end{bmatrix} \prec 0.
 \end{equation*} 
\end{prop}
In IQC-based stability conditions such as Proposition \ref{prop:IQC}, 
it is of prime importance to 
employ a set of multipliers $\mf{\Pi} \subset \mf{\Pi}^{\star}$ that 
can capture the input-output characteristics 
of $\Phi$ as accurately as possible while being numerically tractable. 
As such a set of multipliers, 
this paper focuses on static OZF multipliers for 
slope-restricted and repeated nonlinearities.
% 
%%%%%%%%%%%%%%%%%%%%%%%%%%%%%%%%%%%%%%%%
\subsection{Static O'Shea-Zames-Falb Multipliers}\label{ssec:OZF}
Some specific definitions are necessary to describe static OZF multipliers. 
A matrix $M\in\mb{R}^{m{\times}m}$ is said to be a Z-matrix if 
 $M_{i,j} \le 0$ for all $i \ne j$.
Moreover, $M$ is said to be doubly hyperdominant if
it is a Z-matrix and $M\1 \ge 0$, $\1^TM \ge 0$, where 
$\1\in\mb{R}^m$ stand for the all-ones-vector. 
In addition, $M$ is said to be doubly dominant if 
$|M|_{\mr{d}}\1 \ge 0$, $\1^T|M|_{\mr{d}} \ge 0$. 
In this paper, we denote 
by $\mb{Z}^m,\ \mb{DHD}^m,\ \mb{DD}^m \subset \mb{R}^{m{\times}m}$ the sets of 
Z-matrices, doubly hyperdominant matrices, and doubly dominant matrices, 
respectively. 
Namely, we define 
 \begin{equation*}
   \begin{aligned}
     &
     \mb{DHD}^m := \{
       M\in\mb{Z}^m : M\1 \ge 0, \1^TM \ge 0
       \}, \\
     &
     \mb{DD}^m := \{
       M\in\mb{R}^{m{\times}m} : |M|_{\mr{d}}\1 \ge 0, \1^T|M|_{\mr{d}} \ge 0
       \}.  
   \end{aligned}
 \end{equation*}   
It is obvious that $\mb{DHD}^m \subsetneq \mb{DD}^m$. 
As previously stated, we focus on slope-restricted and repeated nonlinearities
in this paper. 
The definitions of 
the slope-restricted nonlinearity and the sets of repeated nonlinear operators 
are given below.  
\begin{defn} \label{defn:slope-restricted}
Let $\mu \le 0 \le \nu$. 
Then, a nonlinearity $\phi : \mb{R} \rightarrow \mb{R}$ is 
 said to be slope-restricted, in short $\slope{\phi}{\mu}{\nu}$, if 
 $\phi(0) = 0$ and 
 \begin{equation*}
   \mu \le \frac{\phi(p) - \phi(q)}{p - q} \le \nu \ (\forall{p,q}\in\mb{R}, \ p \ne q).
 \end{equation*}
\end{defn}
%%

%% definition of the sets of repeated
\begin{defn} \label{defn:repeated}
For $\mu \le 0 \le \nu$, we define 
 \begin{equation*}
   \begin{aligned}
     &
     \slopePhi{\mu}{\nu} := \{
       \Phi : \Phi = \mr{diag}_m(\phi), \ \slope{\phi}{\mu}{\nu}
       \}, \\
     &
     \oddPhi := \{
       \Phi : \Phi = \mr{diag}_m(\phi), \ \text{$\phi$ is odd}
       \}, 
   \end{aligned}
 \end{equation*}
 where $\mr{diag}_m(\phi) := \mr{diag}(\underbrace{\phi,\ldots,\phi}_m)$.
\end{defn}
Under these definitions, 
the main result of (\cite{Fetzer_IFAC2017}) 
on static OZF multipliers (\cite{O'Shea_IEEE1967,Zames_SIAM1968}) 
for slope-restricted and repeated nonlinearities can be summarized by 
the next lemma.
%

%% OZF basic result
\begin{lem} \label{lem:OZF basic} (\cite{Fetzer_IFAC2017}) \\
Let $\mu \le 0 \le \nu$. 
Then, for all $\Phi\in\slopePhi{\mu}{\nu}$ and $M\in\mb{DHD}^m$, 
 we have
 \begin{equation*}
   \begin{aligned}
     &
     (\ast)^T\ml{M}(M,\mu,\nu)\begin{bmatrix}
       \zeta \\
       \Phi(\zeta)
       \end{bmatrix} \ge 0 \quad \forall{\zeta}\in\mb{R}^m, \\
     &
     \ml{M}(M,\mu,\nu) := (\ast)^T\begin{bmatrix}
       0_m & M \\
       M^T & 0_m
       \end{bmatrix} \begin{bmatrix}
         \nu I_m & -I_m \\
         -\mu I_m & I_m
         \end{bmatrix}.  
   \end{aligned}
 \end{equation*}
Moreover, this result holds for 
 all $\Phi\in\slopePhi{\mu}{\nu}\cap\oddPhi$ and $M\in\mb{DD}^m$. 
\end{lem}
On the basis of this lemma, 
 we define
 \begin{equation*}
   \begin{aligned}
     &
     \mf{\Pi}_{m,\mb{DHD}} := \{
       \Pi\in\mb{S}^{2m} : \Pi = \ml{M}(M,\mu,\nu), \ M\in\mb{DHD}^m
       \}, \\
     & 
     \mf{\Pi}_{m,\mb{DD}} := \{
       \Pi\in\mb{S}^{2m} : \Pi = \ml{M}(M,\mu,\nu), \ M\in\mb{DD}^m
       \}. 
   \end{aligned}
 \end{equation*}
Here, $\mf{\Pi}_{m,\mb{DHD}}$ is used if $\Phi\in\slopePhi{\mu}{\nu}$ and
 $\mf{\Pi}_{m,\mb{DD}}$ is used if $\Phi\in\slopePhi{\mu}{\nu}\cap\oddPhi$.
%   
%%%%%%%%%%%%%%%%%%%%%%%%%%%%%%%%%%%%%%%%%%%%%%%%%%%%%%%%%%%%%%%%%%%%%%%%%
%%%%%%%%%%%%%%%%%%%%%%%%%%%%%%%%%%%%%%%%%%%%%%%%%%%%%%%%%%%%%%%%%%%%%%%%%%
\section{Absolute Stability Analysis for $\Phi\in\slopePhi{\mu}{\nu}$}
%%%%%%%%%%%%%%%%%%%%%%%%%%%%%%%%%%%%%%%%
\subsection{LMIs Ensuring Absolute Stability and Their Dual}
The next lemma readily follows from 
\textit{Proposition \ref{prop:IQC}} and \textit{Lemma \ref{lem:OZF basic}}.
%

%% DHD basic result
\begin{lem} \label{lem:DHD basic}
Let $\mu \le 0 \le \nu$. 
Then, the feedback system $\Sigma$ consisting of 
 \eqref{eq:LTI} and \eqref{eq:nonlinearity} is 
 absolutely stable for $\Phi\in\slopePhi{\mu}{\nu}$ if 
 there exist $P\in\mb{S}_{++}^n$ and $M\in\mb{DHD}^m$ 
 such that
 \begin{equation}
   \primalfirst + \primalsecond{M}{\mu}{\nu} \prec 0.
 \label{eq:pLMI_DHD}
 \end{equation}
\end{lem}
If LMI \eqref{eq:pLMI_DHD} is numerically feasible, 
 we can readily conclude from \textit{Lemma \ref{lem:DHD basic}} that 
 the feedback system $\Sigma$ is absolutely stable for $\Phi\in\slopePhi{\mu}{\nu}$. 
However, since this LMI condition is generally a sufficient condition, 
we cannot conclude anything about the absolute stability of $\Sigma$ 
if the LMI turns out to be numerically infeasible.   
To address this issue, 
we consider the dual of LMI \eqref{eq:pLMI_DHD} that is feasible if and only if 
the (primal) LMI \eqref{eq:pLMI_DHD} is infeasible (\cite{Scherer_EJC2006}). 
In deriving the dual LMI, 
we note that the primal LMI \eqref{eq:pLMI_DHD} can be rewritten equivalently, 
as follows: 
 \\
 \textbf{Primal LMI} (For $\Phi\in\slopePhi{\mu}{\nu}$) % primal for slope[\mu,\nu]
 \\    
 Find $P\in\mb{S}_{++}^n$ and $M\in\mb{R}^{m{\times}m}$ 
  such that
  \begin{equation}
    \begin{aligned}
      &
      \primalfirst + \primalsecond{M}{\mu}{\nu} \prec 0, \\
      &
      M\in\mb{Z}^m,\ M\1 \ge 0,\ \1^TM \ge 0.
    \end{aligned}
  \label{eq:ppLMI_DHD}  
  \end{equation}
To consider the dual of LMI \eqref{eq:ppLMI_DHD},
 we define
 \begin{equation*}
   \mb{Z}_0^m := \{
     X\in\mb{Z}^m : X_{i,i} = 0 \ (i = 1,\ldots,m)
     \}.
 \end{equation*}
Then, for the Lagrange dual variables 
 $H\in\mb{S}_+^{n+m}$, $f,g\in\mb{R}_+^m$ and $X\in\mb{Z}_0^m$, 
 the Lagrangian can be defined as
 \begin{equation*}
 \scalebox{0.85}
 {$
   \begin{aligned}
     &
     \ml{L}(P,M,H,f,g,X) \\
     &
     := \trace \left(
       \left(
         \primalfirst + \primalsecond{M}{\mu}{\nu}
         \right) H 
         \right) \\
     &
     \quad - 2f^TM\1 - 2\1^TMg - 2\trace(MX) \\
     &
     = \trace(P(\dualfirst)) \\
     &
     \quad + 2\trace(M(Y - \1f^T - g\1^T - X))    
   \end{aligned}
 $}  
 \end{equation*}
 where 
 \begin{equation*}
   \begin{aligned}
     H := \begin{bmatrix}
       H_{11} & H_{12} \\
       H_{12}^T & H_{22}
       \end{bmatrix}, \ H_{11}\in\mb{S}_+^n, \ H_{22}\in\mb{S}_+^m, \\
     Y := \matrixY.  
   \end{aligned}
 \end{equation*}
For $\ml{L}(P,M,H,f,g,X) \ge 0$ to hold for 
 any $P\in\mb{S}_{++}^n$ and $M\in\mb{R}^{m{\times}m}$, 
 we can select the solution
 \begin{equation*}
   \dualfirst \succeq 0, \ Y = \1f^T + g\1^T + X.
 \end{equation*} 
We thus arrive at the dual LMI given below.  
\\
\textbf{Dual LMI} (For $\Phi\in\slopePhi{\mu}{\nu}$)
\\
Find $H\in\mb{S}_+^{n+m}$, $f,g\in\mb{R}_+^m$, $X\in\mb{Z}_0^m$, 
 not all zeros, 
 such that
 \begin{equation}
   \begin{aligned}
     &
     \dualfirst \succeq 0, \\
     &
     \matrixY = \1f^T + g\1^T + X.
   \end{aligned}
 \label{eq:dLMI_DHD}  
 \end{equation} 
 %
%
%%%%%%%%%%%%%%%%%%%%%%%%%%%%%%%%%%%%%%%%
\subsection{Detecting Destabilizing Nonlinearity by Dual LMIs}
In this subsection, we consider the dual LMI \eqref{eq:dLMI_DHD}. 
Due to some technical reasons, 
in the following subsections, 
we assume $\|D\| < 1$ and restrict our attention to the case $\mu = 0$ and $\nu = 1$, which yields $\|\Phi\| \le 1$ for $\Phi\in\slopePhi{0}{1}$. 
Then, we can readily ensure the well-posedness of the system $\Sigma$. 
Moreover, by focusing on the upper-left $n{\times}n$ block of \eqref{eq:ppLMI_DHD}, 
we have 
 \begin{equation*}
   -P + A^TPA - \mu\nu C^T(M + M^T)C \prec 0.
 \end{equation*}
It follows that \eqref{eq:ppLMI_DHD} requires $-P + A^TPA \prec 0$ if
 $\mu = 0$ and $\nu = 1$. 
Furthermore, since $A\in\mb{R}^{n{\times}n}$ is assumed to be Schur stable, 
 we see that $P\in\mb{S}_{++}^n$ is automatically satisfied if 
 we simply require $P\in\mb{S}^n$. 
Therefore, 
 the primal LMI \eqref{eq:ppLMI_DHD} and 
 the dual LMI \eqref{eq:dLMI_DHD} reduce respectively to:
 \\
 \textbf{Primal LMI} (For $\Phi\in\slopePhi{0}{1}$)
 \\
 Find $P\in\mb{S}^n$ and $M\in\mb{R}^{m{\times}m}$ 
  such that
  \begin{equation*}
    \begin{aligned}
      &
      \primalfirst + \primalsecond{M}{0}{1} \prec 0, \\
      &
      M\in\mb{Z}^m,\ M\1 \ge 0,\ \1^TM \ge 0.
    \end{aligned}  
  \end{equation*}
 \\ 
 \textbf{Dual LMI} (For $\Phi\in\slopePhi{0}{1}$)
 \\ 
 Find $H\in\mb{S}_+^{n+m}$, $f,g\in\mb{R}_+^m$, $X\in\mb{Z}_0^m$, 
 not all zeros, 
 such that
 \begin{equation}
   \begin{aligned}
     &
     \dualfirst = 0, \\
     &
     \RematrixY = \1f^T + g\1^T + X.
   \end{aligned}
 \label{eq:RedLMI_DHD}  
 \end{equation} 
We note that, in this dual LMI \eqref{eq:RedLMI_DHD}, 
 the first inequality constraint in \eqref{eq:dLMI_DHD} has been replaced by 
 the equality constraint. 
Regarding this dual LMI \eqref{eq:RedLMI_DHD}, 
 we can obtain the next main result for 
 the slope-restricted nonlinearity case.  
%

%% main result for DHD case
\begin{thm} \label{thm:DHD}
Suppose the dual LMI \eqref{eq:RedLMI_DHD} is feasible and
 has a solution $H\in\mb{S}_+^{n+m}$ of $\rank{H} = 1$ given by
 \begin{equation*}
   H = \begin{bmatrix}
     h_1 \\
     h_2
     \end{bmatrix} \begin{bmatrix}
       h_1 \\
       h_2
       \end{bmatrix}^T, \ h_1\in\mb{R}^n, \ h_2\in\mb{R}^m.
 \end{equation*}
We further suppose $\Pd{(Ah_1 + Bh_2) h_1^T} \ge 0$, which implies that the $i$-th elements of the vectors $Ah_1 + Bh_2$ and $h_1$ have the same sign for each $i=1,\ldots,n$. 
Then, the following assertions hold: 
 \begin{enumerate}
 \renewcommand{\labelenumi}{(\roman{enumi})}
 \item \ 
   $h_1 \ne 0$.
 \item \ 
   Let us define $\z := Ch_1 + Dh_2$ and $\w := h_2$. 
   Then, there exists $\slope{\phiwc}{0}{1}$ such that 
    $\phiwc(\z_i) = \w_i \ (i = 1,\ldots,m)$.  
 \item \ 
   Let us define $\Phiwc := \mr{diag}_m(\phiwc)$. 
   Then, the feedback system $\Sigma$ with 
    the nonlinearity $\Phi = \Phiwc\in\slopePhi{0}{1}$ is unstable. 
   In particular, $h_1$ is a nonzero equilibrium point of $\Sigma$, i.e, $x(k) = h_1 \ (k=1,2,\ldots)$ if $x(0) = h_1$. 
 \end{enumerate} 
\end{thm}
%%

%% proof of theorem (DHD)
\begin{proofof}{\textit{Theorem \ref{thm:DHD}}} 

%=========================================================================
\underline{\textit{Proof of} (i)} : \ 
To prove that $\rank{H} = 1$ can happen,
 we first prove that $H \ne 0$.
To this end, suppose $H = 0$ for proof by contradiction. 
Then, from the second equality constraint in \eqref{eq:RedLMI_DHD}, 
 we see $f = 0$, $g = 0$, and $X = 0$, 
 which contradicts the requirements that $H,f,g,X$ are not all zeros. 
Therefore $H \ne 0$. 
We next suppose $h_1 = 0$ for proof by contradiction. 
Then, $h_2 \ne 0$ since $H \ne 0$ as proved. 
Again from the second equality constraint in \eqref{eq:RedLMI_DHD},
 one has $\trace(h_2h_2^T(D^T - I_m)) \ge 0$ and 
 this implies $h_2^TDh_2 \ge h_2^Th_2$. 
However, this is impossible for $h \ne 0$ since
 we assumed $\|D\| < 1$. 
Therefore $h_1 \ne 0$.  
%
%=========================================================================

%=========================================================================
\underline{\textit{Proof of} (ii)} : \ 
We first prove that $\phiwc : \mb{R} \rightarrow \mb{R}$ satisfying
 $\phiwc(\z_i) = \w_i \ (i = 1,\ldots,m)$ is well-defined. 
To this end, it suffices to prove that
 if $\z_i = \z_j$ then $\w_i = \w_j \ (i \ne j)$. 
To prove this assertion,
 we note that 
 the second equality constraint in \eqref{eq:RedLMI_DHD} can be rewritten 
 equivalently as
 \begin{equation}
   \w(\z - \w)^T = \1f^T + g\1^T + X.
 \label{eq:proof_DHD_2.1}  
 \end{equation} 
Then, if $\z_i = \z_j \ (i \ne j)$, we obtain from \eqref{eq:proof_DHD_2.1} that
 \begin{equation*}
   \begin{aligned}
     -(\w_i - \w_j)^2 
     & = (e_i -e_j)^T\w(\z - \w)^T(e_i -e_j) \\
     & = (e_i -e_j)^T(\1f^T + g\1^T + X)(e_i -e_j) \ge 0
   \end{aligned}
 \end{equation*} 
 where $e_i\in\mb{R}^m \ (i = 1,\ldots,m)$ are the standard basis. 
This clearly shows $\w_i = \w_j \ (i \ne j)$. 
We finally prove the existence of
 $\slope{\phiwc}{0}{1}$ such that $\phiwc(\z_i) = \w_i \ (i = 1,\ldots,m)$. 
To this end, it is first needed to prove that
 if $\z_i = 0$ then $\w_i = 0$ from \textit{Definition \ref{defn:slope-restricted}}. 
This can be verified from \eqref{eq:proof_DHD_2.1} since
 if $\z_i = 0$ then we obtain
 \begin{equation*}
   \begin{aligned}
     -w_i^{{\ast}2} 
     & = e_i^T\w(\z - \w)^Te_i \\
     & = e_i^T(\1f^T + g\1^T + X)e_i \ge 0.
   \end{aligned}
 \end{equation*} 
This clearly shows $\w_i = 0$. 
Next, we verify that the inequality in \textit{Definition \ref{defn:slope-restricted}} holds. 
Again from \eqref{eq:proof_DHD_2.1},
 we obtain
 \begin{equation*}
   \begin{aligned}
     &
     (\w_i - \w_j)((\z_i - \z_j) - (\w_i - \w_j)) \\
     & 
     = (e_i - e_j)^T\w(\z - \w)^T(e_i - e_j) \\
     &
     = (e_i - e_j)^T(\1f^T + g\1^T + X)(e_i - e_j) \ge 0.
   \end{aligned}  
 \end{equation*}
Therefore, if $\z_i \ne \z_j \ (i \ne j)$,
 we have
 \begin{equation}
   \frac{\w_i - \w_j}{\z_i - \z_j} \left(
     1 - \frac{\w_i - \w_j}{\z_i - \z_j}
     \right) \ge 0.
 \label{eq:proof_DHD_2.2}    
 \end{equation}
In addition, from \eqref{eq:proof_DHD_2.1}, 
 we obtain
 \begin{equation*}
   \begin{aligned}
     \w_i(\z_i - \w_i) 
     & = e_i^T\w(\z - \w)^Te_i \\
     & = e_i^T(\1f^T + g\1^T + X)e_i \ge 0.
   \end{aligned}
 \end{equation*} 
Therefore, if $\z_i \ne 0 \ (i = 1,\ldots,m)$,
 we have
 \begin{equation}
   \frac{\w_i}{\z_i} \left(
     1 - \frac{\w_i}{\z_i}
     \right) \ge 0.
 \label{eq:proof_DHD_2.3}    
 \end{equation} 
From $\phiwc(0) = 0$, \eqref{eq:proof_DHD_2.2}, and \eqref{eq:proof_DHD_2.3},
 we can conclude that 
 $\slope{\phiwc}{0}{1}$ such that $\phiwc(\z_i) = \w_i \ (i = 1,\ldots,m)$ exists.  
% 
%========================================================================= 
 
%=========================================================================
\underline{\textit{Proof of} (iii)} : \ 
To validate the assertion (iii),
 it suffices to show that 
 $x(k) = h_1$, $z(k) = \z$, and $w(k) = \w \ (k=0,1,2,\cdots)$ satisfy 
 \eqref{eq:LTI} and \eqref{eq:nonlinearity} in 
 the case $\Phi = \Phiwc$ and $x(0) = h_1$.
From the first equality constraint in \eqref{eq:RedLMI_DHD},
 we have
 \begin{equation*}
   -h_1h_1^T + (Ah_1 + Bh_2)(Ah_1 + Bh_2)^T = 0.
 \end{equation*} 
Since $h_1 \ne 0$ as proved,
 we see from (\cite{green}) 
 that $Ah_1 + Bh_2 = h_1$ or $Ah_1 + Bh_2 = -h_1$. 
Furthermore, from the assumption that $\Pd{(Ah_1 + Bh_2)h_1^T} \ge 0$,
 we have $Ah_1 + Bh_2 = h_1$. 
Therefore, 
 we can obtain
 \begin{equation*}
   \begin{aligned}
     &
     x(k+1) = h_1 = Ah_1 + Bh_2 = Ax(k) + B\w, \\
     &
     z(k) = \z = Ch_1 + Dh_2 = Cx(k) + D\w, \\
     &
     w(k) = \w = \Phiwc(\z).
   \end{aligned}
 \end{equation*}    
This completes the proof.  
\end{proofof}
%%

%%%%%%%%%%%%%%%%%%%%%%%%%%%%%%%%%%%%%%%%
\subsection{Concrete Construction of Destabilizing Nonlinearity}\label{ssec:Construction_DHD}
From \textit{Theorem \ref{thm:DHD}},
 we see that any $\slope{\phi}{0}{1}$ with $\phi(\z_i) = \w_i \ (i = 1,\ldots,m)$ is 
 a destabilizing nonlinearity. 
One of such destabilizing nonlinear (piecewise linear) operators can be constructed by 
 following the next procedure: 
 \begin{enumerate}
 \renewcommand{\labelenumi}{\arabic{enumi})}
 \item \ 
   Define the set $\ml{Z}_0 := \{ 0,\z_1,\z_2,\ldots,\z_m \}$ (by
    leaving only one if they have duplicates) and 
    compose the series $\bar{z}_1,\bar{z}_2,\ldots,\bar{z}_l$, $l = |\ml{Z}_0|$, 
    such that each $\bar{z}_i \ (1 \le i \le l)$ is 
    the $i$-th smallest value of $\ml{Z}_0$. 
   Similarly, define the series $\bar{w}_1,\bar{w}_2,\ldots,\bar{w}_l$ for 
   the set $\{ 0,\w_1,\w_2,\ldots,\w_m \}$.
 \item \ 
   Define $\phiwc$ as follows: 
    \begin{equation*}
      \begin{aligned}
        &
        \phiwc(z) \\
        &
        = \begin{cases}
          \bar{w}_1 & z < \bar{z}_1, \\
          \frac{\bar{w}_{i+1} - \bar{w}_i}{\bar{z}_{i+1} - \bar{z}_i} (z - \bar{z}_i) + \bar{w}_i
           & \vspace*{-2mm} \bar{z}_i \le z \le \bar{z}_{i+1} \\
           & (i = 1,\ldots,l-1), \\
          \bar{w}_l & \bar{z}_l \le z. 
          \end{cases}
      \end{aligned}
    \end{equation*}    
 \end{enumerate}
%
%%%%%%%%%%%%%%%%%%%%%%%%%%%%%%%%%%%%%%%%
\subsection{Numerical Examples}   
\begin{figure}[t]
  \centering
  \includegraphics[scale = 0.6]{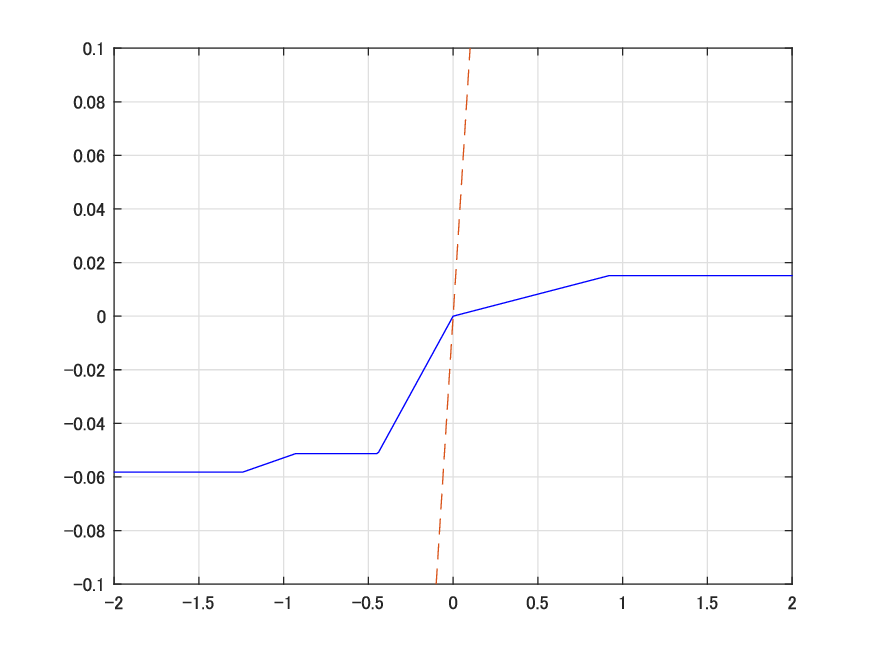}
  \vspace*{-10mm}
  \caption{The input-output map of the detected $\phiwc$}
  \label{fig:in-out_DHD}
  \centering\hspace*{-6mm}
  \includegraphics[scale = 0.65]{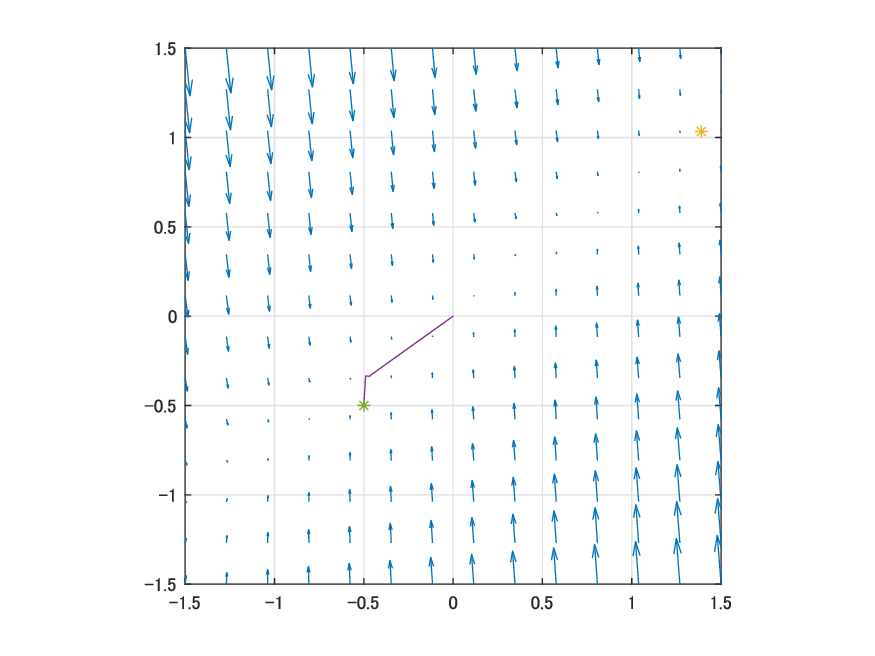}
  \vspace*{-10mm}
  \caption{State trajectories with $x(0) = h_1$ and $x(0) = [-0.5\ -0.5]^T$}
  \label{fig:trajectory_DHD}
\end{figure}

In this section,
 we demonstrate the effectiveness of the result in 
 \textit{Theorem \ref{thm:DHD}} by a numerical example. 
Let us consider the case where 
 the coefficient matrices in \eqref{eq:LTI} are given by
 \begin{equation*}
   \begin{aligned}
     A = \begin{bmatrix*}[r] % matrix A
       0.88 & 0.06 \\
       0.73 & -0.05
       \end{bmatrix*}, \ & B = \begin{bmatrix*}[r] % matrix B
         -0.49 & -0.65 & -0.75 & 0.22 \\
         0.03 & -0.86 & -0.53 & -0.32
         \end{bmatrix*}, \\
     C = \begin{bmatrix*}[r] % matrix C
       -0.27 & -0.55 \\
       -0.78 & -0.09 \\
       -0.23 & -0.09 \\
       0.46 & 0.27
       \end{bmatrix*}, \ & D = \begin{bmatrix*}[r] % matrix D
         0.39 & -0.63 & 0.03 & 0.06 \\
         0.11 & 0.72 & 0.25 & -0.28 \\
         -0.20 & 0.60 & -0.14 & -0.91 \\
         -0.20 & 0.72 & -0.68 & -0.04
         \end{bmatrix*}.    
   \end{aligned}
 \end{equation*} 
For this system,
 the dual LMI \eqref{eq:RedLMI_DHD} turns out to be feasible, 
 and the resulting dual solution $H$ is numerically verified to 
 be $\rank{H} = 1$. 
The full-rank factorization of $H$ as well as
 $\z,\w\in\mb{R}^4$ in \textit{Theorem \ref{thm:DHD}} are obtained as
 \begin{equation*}
 \scalebox{0.7}
 {$
   h_1 = \begin{bmatrix*}[r]
     1.3898 \\
     1.0337
     \end{bmatrix*}, \ h_2 = \begin{bmatrix*}[r]
       -0.0513 \\
       -0.0582 \\
       -0.0513 \\
       0.0151
       \end{bmatrix*}, \ \z = \begin{bmatrix*}[r]
         -0.9277 \\
         -1.2417 \\
         -0.4440 \\
         0.9210
         \end{bmatrix*}, \ \w = \begin{bmatrix*}[r]
           -0.0513 \\
           -0.0582 \\
           -0.0513 \\
           0.0151
           \end{bmatrix*} (= h_2).
 $}          
 \end{equation*}  
It is obvious that $h_1 \ne 0$ (the assertion (i) of \textit{Theorem \ref{thm:DHD}}). 
In Fig.~\ref{fig:in-out_DHD},
 the solid line represents the input-output characteristics of 
 the destabilizing nonlinear (piecewise linear) operator $\phiwc$ constructed by 
 following the procedure in the preceding subsection, 
 while the dashed line represents a line of slope one. 
From Fig.~\ref{fig:in-out_DHD}, 
 we can readily see $\slope{\phiwc}{0}{1}$ (the assertion (ii)). 
Fig.~\ref{fig:trajectory_DHD} shows the vector field of 
 the feedback system $\Sigma$ with the nonlinearity $\Phi = \Phiwc$, 
 together with the state trajectories from 
 initial states $x(0) = h_1$ and $x(0) = [-0.5\ -0.5]^T$. 
The state trajectory from
 the initial state $x(0) = [-0.5\ -0.5]^T$ converges to the origin.   
However, as proved in \textit{Theorem \ref{thm:DHD}},
 the state trajectory from 
 the initial state $x(0) = h_1$ does not evolve and 
 $x(0) = h_1$ is confirmed to be an equilibrium point of 
 the feedback system $\Sigma$ (the assertion (iii)).
%
%%%%%%%%%%%%%%%%%%%%%%%%%%%%%%%%%%%%%%%%%%%%%%%%%%%%%%%%%%%%%%%%%%%%%%%%%
%%%%%%%%%%%%%%%%%%%%%%%%%%%%%%%%%%%%%%%%%%%%%%%%%%%%%%%%%%%%%%%%%%%%%%%%%
\section{Absolute Stability Analysis for $\Phi\in\slopePhi{\mu}{\nu}\cap\oddPhi$}
\label{sec:slope-odd}
%%%%%%%%%%%%%%%%%%%%%%%%%%%%%%%%%%%%%%%%
\subsection{LMIs Ensuring Absolute Stability and Their Dual}
The next lemma readily follows from 
\textit{Proposition \ref{prop:IQC}} and \textit{Lemma \ref{lem:OZF basic}}.

%% DD basic result
\begin{lem} \label{lem:DD basic}
Let $\mu \le 0 \le \nu$. 
Then, the feedback system $\Sigma$ consisting of 
 \eqref{eq:LTI} and \eqref{eq:nonlinearity} is 
 absolutely stable for $\Phi\in\slopePhi{\mu}{\nu}\cap\oddPhi$ if 
 there exist $P\in\mb{S}_{++}^n$ and $M\in\mb{DD}^m$ 
 such that
 \begin{equation}
   \primalfirst + \primalsecond{M}{\mu}{\nu} \prec 0.
 \label{eq:pLMI_DD}
 \end{equation}
\end{lem}
For the same reasons as the slope-restricted and repeated nonlinearities case,
we derive the dual of LMI \eqref{eq:pLMI_DD}. 
To this end, we first note that LMI \eqref{eq:pLMI_DD} can be rewritten equivalently, 
as follows: 
 \\
 \textbf{Primal LMI} (For $\Phi\in\slopePhi{\mu}{\nu}\cap\oddPhi$) 
 \\ 
 Find $P\in\mb{S}_{++}^n$, $\Md\in\mb{D}^m$, and $\Mod,\oMod\in\mb{OD}^m$
  such that
  \begin{equation}
  \scalebox{0.9}
  {$
    \begin{aligned}
      &
      \primalfirst + \primalsecond{\Md + \Mod}{\mu}{\nu} \prec 0, \\
      &
      (\Md - \oMod)\1 \ge 0, \ \1^T(\Md - \oMod) \ge 0, \\
      &
      \oMod - \Mod \ge 0, \ \oMod + \Mod \ge 0.
    \end{aligned}
  $}  
  \label{eq:ppLMI_DD}  
  \end{equation}
Then, for 
 the Lagrange dual variables $H\in\mb{S}_+^{n+m}$, $f,g\in\mb{R}_+^m$, 
 and $X,Z\in\mb{Z}_0^m$, 
 the Lagrangian can be defined as
 \begin{equation*}
 \scalebox{0.75}
 {$
   \begin{aligned}
     &
     \ml{L}(P,\Md,\Mod,\oMod,H,f,g,X,Z) \\
     &
     := \trace \left(
       \left(
         \primalfirst + \primalsecond{\Md + \Mod}{\mu}{\nu}
         \right) H
         \right) \\
     &
     \quad -2f^T(\Md - \oMod)\1 - 2\1^T(\Md - \oMod)g \\
     &
     \quad + 2\trace((\oMod - \Mod)X) + 2\trace((\Mod + \oMod)Z) \\
     &
     = \trace(P(\dualfirst)) \\
     &
     \quad + 2\trace(\Md(Y - \1f^T - g\1^T)) \\
     &
     \quad + 2\trace(\Mod(Y - X + Z)) \\
     &
     \quad + 2\trace(\oMod(\1f^T + g\1^T + X + Z))    
   \end{aligned}
 $}  
 \end{equation*}
 where
 \begin{equation*}
   \begin{aligned}
     &
     H := \begin{bmatrix}
       H_{11} & H_{12} \\
       H_{12}^T & H_{22}
       \end{bmatrix}, \ H_{11}\in\mb{S}_+^n, \ H_{22}\in\mb{S}_+^m, \\
     &
     Y := \matrixY.  
   \end{aligned}
 \end{equation*}
For $\ml{L}(P,\Md,\Mod,\oMod,H,f,g,X,Z) \ge 0$ to hold for
 any $P\in\mb{S}_{++}^n$, $\Md\in\mb{D}^m$, and $\Mod,\oMod\in\mb{OD}^m$, 
 we see
 \begin{equation*}
   \begin{aligned}
     &
     \dualfirst \succeq 0, \\
     &
     \Pd{Y} = \Pd{\1f^T + g\1^T}, \\
     &
     \Pod{Y} = \Pod{X - Z}, \\
     &
     \Pod{X + Z} = -\Pod{\1f^T + g\1^T}.
   \end{aligned}
 \end{equation*}
We thus arrive at the dual LMI given below.  
\\
\textbf{Dual LMI} (For $\Phi\in\slopePhi{\mu}{\nu}\cap\oddPhi$)
\\ 
Find $H\in\mb{S}_+^{n+m}$, $f,g\in\mb{R}_+^m$, $X,Z\in\mb{Z}_0^m$,
 not all zeros, such that
 \begin{equation}
   \begin{aligned}
     &
     \dualfirst \succeq 0, \\
     &
     \Pd{\matrixY} \\
     &
     \quad = \Pd{\1f^T + g\1^T}, \\
     &
     \Pod{\matrixY} \\
     &
     \quad = \Pod{X - Z}, \\
     &
     \Pod{X + Z} = -\Pod{\1f^T + g\1^T}.
   \end{aligned}
 \label{eq:dLMI_DD}  
 \end{equation}
 %
%
%%%%%%%%%%%%%%%%%%%%%%%%%%%%%%%%%%%%%%%%
\subsection{Detecting Destabilizing Nonlinearity by Dual LMIs}
Suppose $\mu = 0$, $\nu = 1$, and $\|D\| < 1$ 
for the same reasons as the slope-restricted and repeated nonlinearities case. 
Then, 
 the primal LMI \eqref{eq:ppLMI_DD} and 
 the dual LMI \eqref{eq:dLMI_DD} reduce respectively to:
 \\
 \textbf{Primal LMI} (For $\Phi\in\slopePhi{0}{1}\cap\oddPhi$) 
 \\
 Find $P\in\mb{S}^n$, $\Md\in\mb{D}^m$, and $\Mod,\oMod\in\mb{OD}^m$
  such that
  \begin{equation*}
  \scalebox{0.9}
  {$
    \begin{aligned}
      &
      \primalfirst + \primalsecond{\Md + \Mod}{0}{1} \prec 0, \\
      &
      (\Md - \oMod)\1 \ge 0, \ \1^T(\Md - \oMod) \ge 0, \\
      &
      \oMod - \Mod \ge 0, \ \oMod + \Mod \ge 0.
    \end{aligned}
  $}  
  \end{equation*}
 \\
 \textbf{Dual LMI} (For $\Phi\in\slopePhi{0}{1}\cap\oddPhi$)
 \\ 
 Find $H\in\mb{S}_+^{n+m}$, $f,g\in\mb{R}_+^m$, $X,Z\in\mb{Z}_0^m$,
  not all zeros, such that
  \begin{equation}
    \begin{aligned}
      &
      \dualfirst = 0, \\
      &
      \Pd{\RematrixY} \\
      &
      \quad = \Pd{\1f^T + g\1^T}, \\
      &
      \Pod{\RematrixY} \\
      &
      \quad = \Pod{X - Z}, \\
      &
      \Pod{X + Z} = -\Pod{\1f^T + g\1^T}.
    \end{aligned}
  \label{eq:RedLMI_DD}  
  \end{equation}
Regarding this dual LMI \eqref{eq:RedLMI_DD},
 we can obtain the next main result for 
 the slope-restricted and odd nonlinearity case.
%

%% main result for DD case
\begin{thm} \label{thm:DD}
Suppose the dual LMI \eqref{eq:RedLMI_DD} is feasible and
 has a solution $H\in\mb{S}_+^{n+m}$ of $\rank{H} = 1$ given by
 \begin{equation*}
   H = \begin{bmatrix}
     h_1 \\
     h_2
     \end{bmatrix} \begin{bmatrix}
       h_1 \\
       h_2
       \end{bmatrix}^T, \ h_1\in\mb{R}^n, \ h_2\in\mb{R}^m.
 \end{equation*}
We further suppose $\Pd{(Ah_1 + Bh_2) h_1^T} \ge 0$, which implies that the $i$-th elements of the vectors $Ah_1 + Bh_2$ and $h_1$ have the same sign for each $i=1,\ldots,n$. 
Then, the following assertions hold: 
 \begin{enumerate}
 \renewcommand{\labelenumi}{(\roman{enumi})}
 \item \ 
   $h_1 \ne 0$.
 \item \ 
   Let us define $\z := Ch_1 + Dh_2$ and $\w := h_2$. 
   Then, there exists an odd nonlinearity $\slope{\phiwc}{0}{1}$ such that 
    $\phiwc(\z_i) = \w_i$ and $\phiwc(-\z_i) = -\w_i \ (i = 1,\ldots,m)$.  
 \item \ 
   Let us define $\Phiwc := \mr{diag}_m(\phiwc)$. 
   Then, the feedback system $\Sigma$ with 
    the nonlinearity $\Phi = \Phiwc\in\slopePhi{0}{1}\cap\oddPhi$ is unstable. 
   In particular, $h_1$ is a nonzero equilibrium point of $\Sigma$, i.e, $x(k) = h_1 \ (k=1,2,\ldots)$ if $x(0) = h_1$. 
 \end{enumerate} 
\end{thm}
%%  

%% proof of theorem (DD)
\begin{proofof}{\textit{Theorem \ref{thm:DD}}}

%========================================================================
\underline{\textit{Proof of} (i)} : \ 
To prove that $\rank{H} = 1$ can happen,
 we first prove that $H \ne 0$. 
To this end, suppose $H = 0$ for proof by contradiction. 
Then, from the second equality constraint in \eqref{eq:RedLMI_DD},
 we see $f = 0$ and $g = 0$, 
 and hence from the fourth equality constraint in \eqref{eq:RedLMI_DD}, 
 we see $X = 0$ and $Z = 0$, 
 which contradicts the requirements that 
 $H,f,g,X,Z$ are not all zeros. 
Therefore $H \ne 0$. 
We next suppose $h_1 = 0$ for proof by contradiction. 
Then, $h_2 \ne 0$ since $H \ne 0$ as proved. 
We see $\trace(h_2h_2^T(D^T - I_m)) \ge 0$ from 
 the second equality constraint in \eqref{eq:RedLMI_DD}. 
This implies $h_2^TDh_2 \ge h_2^Th_2$,
 which does not hold for $h_2 \ne 0$ since we assumed $\|D\| < 1$. 
Therefore $h_1 \ne 0$.
%
%========================================================================

%========================================================================
\underline{\textit{Proof of} (ii)} : \ 
We first prove that the odd function $\phiwc : \mb{R} \rightarrow \mb{R}$ satisfying
 $\phiwc(\z_i) = \w_i$ and $\phiwc(-\z_i) = -\w_i \ (i = 1,\ldots,m)$ is well-defined. 
To this end, it suffices to prove that
 (iia) \ if $\z_i = 0$ then $\w_i = 0$, 
 (iib) \ if $\z_i = \z_j$ then $\w_i = \w_j \ (i \ne j)$, 
 and (iic) \ if $\z_i = -\z_j$ then $\w_i = -\w_j \ (i \ne j)$.
To prove these assertions,
 we note that the second and third equality constraint in \eqref{eq:RedLMI_DD} 
 can be rewritten equivalently and respectively as
 \begin{align}
   &
   \Pd{\w(\z - \w)^T} = \Pd{\1f^T + g\1^T}, \label{eq:proof_DD_2.1} \\
   &
   \Pod{\w(\z - \w)^T} = \Pod{X - Z}. \label{eq:proof_DD_2.2}
 \end{align}
To prove (iia), suppose $\z_i = 0$. 
Then, we obtain from \eqref{eq:proof_DD_2.1} that
 \begin{equation*}
   -w_i^{{\ast}2} = e_i^T\w(\z - \w)^Te_i = e_i^T(\1f^T + g\1^T)e_i \ge 0.
 \end{equation*}
This clearly shows $\w_i = 0$. 
To prove (iib), suppose $\z_i = \z_j \ (i \ne j)$. 
Then,
 from \eqref{eq:proof_DD_2.1}, \eqref{eq:proof_DD_2.2}, and 
 the fourth equality constraint in \eqref{eq:RedLMI_DD}, 
 we have
 \begin{equation*}
   \begin{aligned}
     & \hspace*{-5mm}
     -(\w_i - \w_j)^2 \\
     &
     = (e_i - e_j)^T\w(\z - \w)^T(e_i - e_j) \\
     &
     = e_i^T(\1f^T + g\1^T)e_i + e_j(\1f^T + g\1^T)e_j \\
     &
     \quad -e_j^T(X - Z)e_i - e_i^T(X - Z)e_j \\
     &
     \ge e_i^T(\1f^T + g\1^T)e_i + e_j(\1f^T + g\1^T)e_j \\
     &
     \quad +e_j^T(X + Z)e_i + e_i^T(X + Z)e_j \\
     &
     = e_i^T(\1f^T + g\1^T)e_i + e_j(\1f^T + g\1^T)e_j \\
     &
     \quad -e_j^T(\1f^T + g\1^T)e_i - e_i^T(\1f^T + g\1^T)e_j \\
     &
     = f_i + g_i + f_j + g_j - (f_i + g_j) - (f_j + g_i) = 0.
   \end{aligned}
 \end{equation*}
This clearly shows $\w_i = \w_j \ (i \ne j)$. 
To prove (iic), suppose $\z_i = -\z_j \ (i \ne j)$. 
Then, similarly to the proof of (iib),
 we have
 \begin{equation*}
   \begin{aligned}
     & \hspace*{-5mm}
     -(\w_i + \w_j)^2 \\
     &
     = (e_i + e_j)^T\w(\z - \w)^T(e_i + e_j) \\
     &
     = e_i^T(\1f^T + g\1^T)e_i + e_j^T(\1f^T + g\1^T)e_j \\
     &
     \quad + e_j^T(X - Z)e_i + e_i^T(X - Z)e_j \\
     &
     \ge e_i^T(\1f^T + g\1^T)e_i + e_j^T(\1f^T + g\1^T)e_j \\
     &
     \quad + e_j^T(X + Z)e_i + e_i^T(X + Z)e_j = 0
   \end{aligned}
 \end{equation*} 
This clearly shows $\w_i = -\w_j \ (i \ne j)$.  
We finally prove the existence of odd $\slope{\phiwc}{0}{1}$ such that
 $\phiwc(\z_i) = \w_i$ and $\phiwc(-\z_i) = \w_i \ (i = 1,\ldots,m)$. 
From \eqref{eq:proof_DD_2.1} and \eqref{eq:proof_DD_2.2},
 we obtain 
 \begin{equation*}
   \begin{aligned}
     &
     (\w_i - \w_j)((\z_i - \z_j) - (\w_i - \w_j)) \\
     &
     = (e_i - e_j)^T\w(\z - \w)^T(e_i - e_j) \ge 0.
   \end{aligned}
 \end{equation*} 
Therefore, if $\z_i \ne \z_j \ (i \ne j)$,
 we have
 \begin{equation}
   \frac{\w_i - \w_j}{\z_i - \z_j} \left(
     1 - \frac{\w_i - \w_j}{\z_i - \z_j}
     \right) \ge 0.
 \label{eq:proof_DD_2.3}
 \end{equation} 
Again from \eqref{eq:proof_DD_2.1} and \eqref{eq:proof_DD_2.2},
 we namely obtain
 \begin{equation*}
   \begin{aligned}
     &
     (\w_i + \w_j)((\z_i + \z_j) - (\w_i + \w_j)) \\
     &
     = (e_i + e_j)^T\w(\z - \w)^T(e_i + e_j) \ge 0.
   \end{aligned}
 \end{equation*} 
 Therefore, if $\z_i \ne -\z_j \ (i \ne j)$,
 we have
 \begin{equation}
   \frac{\w_i - (-\w_j)}{\z_i - (-\z_j)} \left(
     1 - \frac{\w_i - (-\w_j)}{\z_i - (-\z_j)}
     \right) \ge 0.
 \label{eq:proof_DD_2.4}
 \end{equation} 
Moreover, from \eqref{eq:proof_DD_2.1},
 we obtain
 \begin{equation*}
   \w_i(\z_i - \w_i) = e_i^T\w(\z - \w)^Te_i = e_i^T(\1f^T + g\1^T)e_i \ge 0.
 \end{equation*} 
Therefore, if $\z_i \ne 0 \ (i = 1,\ldots,m)$,
 we have
 \begin{equation}
   \frac{\w_i}{\z_i} \left(
     1 - \frac{\w_i}{\z_i}
     \right) \ge 0.
 \label{eq:proof_DD_2.5}    
 \end{equation} 
From $\phiwc(0) = 0$, \eqref{eq:proof_DD_2.3}, \eqref{eq:proof_DD_2.4}, 
 and \eqref{eq:proof_DD_2.5}, 
 we can conclude that an odd $\slope{\phiwc}{0}{1}$ such that 
 $\phiwc(\z_i) = \w_i$ and $\phiwc(-\z_i) = -\w_i \ (i = 1,\ldots,m)$ exists.
%
%========================================================================

%========================================================================
\underline{\textit{Proof of} (iii)} : \ 
Omitted since the proof is exactly the same as
 the proof of (iii) in \textit{Theorem \ref{thm:DHD}}.
\end{proofof}

\begin{rem}
As stated in Section~\ref{sec:intro}, 
the results in Theorems~\ref{thm:DHD} and \ref{thm:DD} can be regarded as
extensions of those obtained for 
continuous-time feedback systems (\cite{Gyotoku_ECC2025})
to discrete-time feedback systems.   
In stark contrast with the
continuous-time case results where only $\rank{H}=1$ is 
required to detect destabilizing nonlinearities  (\cite{Gyotoku_ECC2025}), 
here we are further required to satisfy $\Pd{(Ah_1 + Bh_2) h_1^T} \ge 0$
on the dual solution.  
The system theoretic interpretation of this condition, 
as well as the active use of this condition to enhance the detection of
destabilizing nonlinearities, are currently under investigation.  
\label{rem:novel}
\end{rem}

%%%%%%%%%%%%%%%%%%%%%%%%%%%%%%%%%%%%%%%%
\subsection{Concrete Construction of Destabilizing Nonlinearity}\label{ssec:Construction_DD}
From \textit{Theorem \ref{thm:DD}},
 we see that any odd $\slope{\phi}{0}{1}$ with 
 $\phiwc(\z_i) = \w_i$ and $\phiwc(-\z_i) = -\w_i \ (i = 1,\ldots,m)$ is 
 a destabilizing nonlinearity. 
One of such destabilizing nonlinear (piecewise linear) operators can be constructed by
 following the next procedure:
 \begin{enumerate}
 \renewcommand{\labelenumi}{\arabic{enumi})}
 \item \ 
   Define the set
    \begin{equation*}
      \ml{Z}_0 := \{ 0,\z_1,-\z_1,\z_2,-\z_2,\ldots,\z_m,-\z_m \}
    \end{equation*}
    (by leaving only one if they have duplicates) and 
    compose the series $\bar{z}_1,\bar{z}_2,\ldots,\bar{z}_l, \ l = |\ml{Z}_0|$, 
    such that each $\bar{z}_i \ (1 \le i \le l)$ is 
    the \textit{i}-th smallest value of $\ml{Z}_0$. 
   Similarly, define the series $\bar{w}_1,\bar{w}_2,\ldots,\bar{w}_l$ for
    the set $\{ 0,\w_1,-\w_1,\w_2,-\w_2,\ldots,\w_m,-\w_m \}$. 
 \item \ 
   Define $\phiwc$ as follows:
    \begin{equation*}
      \begin{aligned}
        &
        \phiwc(z) \\
        &
        = \begin{cases}
          \bar{w}_1 & z < \bar{z}_1, \\
          \frac{\bar{w}_{i+1} - \bar{w}_i}{\bar{z}_{i+1} - \bar{z}_i} (z - \bar{z}_i) + \bar{w}_i
           & \bar{z}_i \le z \le \bar{z}_{i+1} \\
           & (i = 1,\ldots,l-1), \\
          \bar{w}_l & \bar{z}_l \le z. 
          \end{cases}
      \end{aligned}
    \end{equation*}    
 \end{enumerate}
%   
%%%%%%%%%%%%%%%%%%%%%%%%%%%%%%%%%%%%%%%%
\subsection{Numerical Examples}
\begin{figure}[t]
  \centering
  \includegraphics[scale = 0.6]{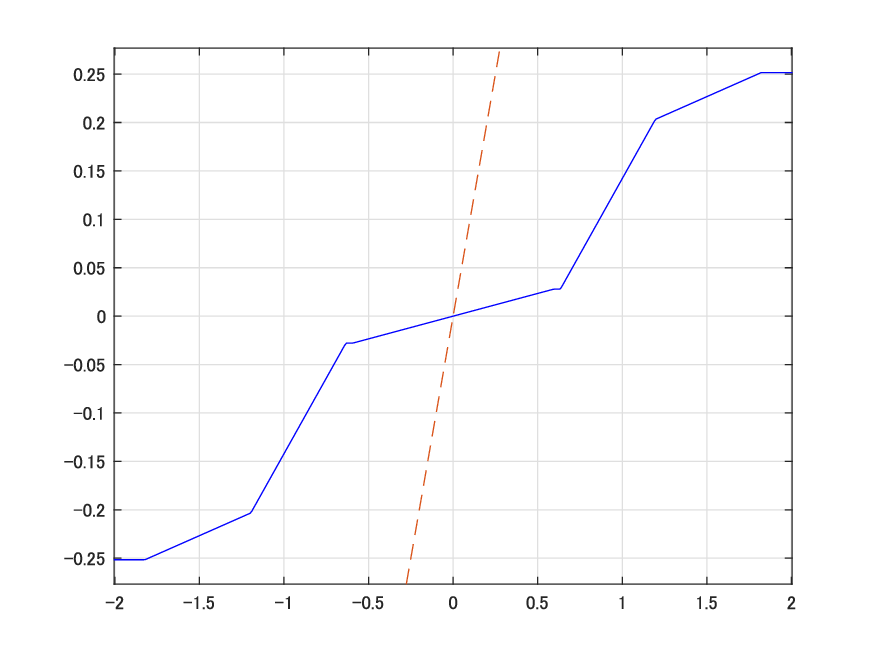}
  \vspace*{-10mm}
  \caption{The input-output map of the detected $\phiwc$}
  \label{fig:in-out_DD}
  \centering\hspace*{-6mm}
  \includegraphics[clip,scale = 0.65]{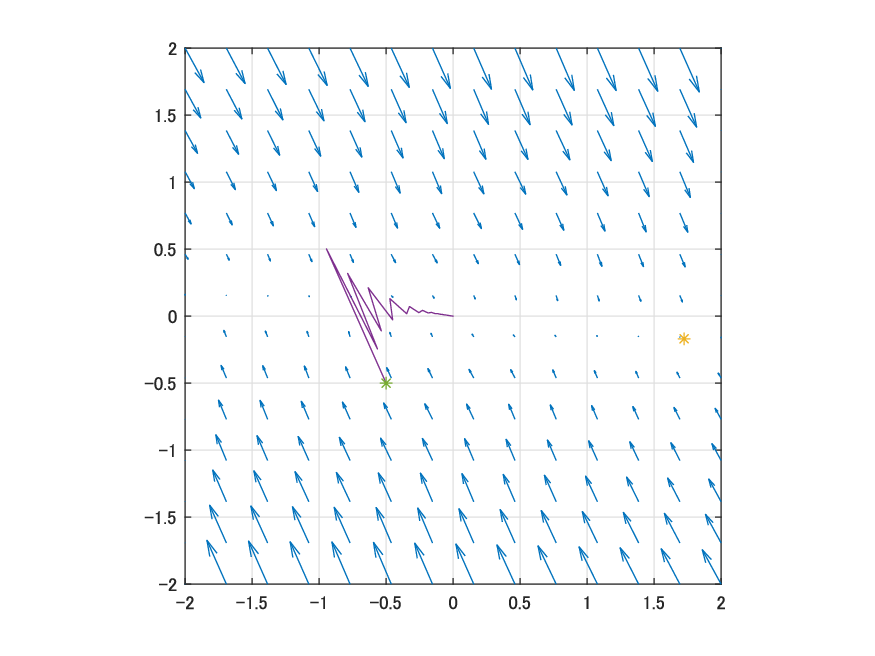}
  \vspace*{-10mm}
  \caption{State trajectories with $x(0) = h_1$ and $x(0) = [-0.5\ -0.5]^T$}
  \label{fig:trajectory_DD}
\end{figure}

In this section,
 we demonstrate the effectiveness of the result in 
 \textit{Theorem \ref{thm:DD}} by a numerical example. 
Let us consider the case where 
 the coefficient matrices in \eqref{eq:LTI} are given by
 \begin{equation*}
   \begin{aligned}
     A = \begin{bmatrix*}[r] % matrix A
       0.95 & 0.95 \\
       -0.28 & -0.72
       \end{bmatrix*}, \ & B = \begin{bmatrix*}[r] % matrix B
         -0.49 & -0.18 & 0.59 & -0.05 \\
         -0.63 & -0.02 & 0.15 & 0.08
         \end{bmatrix*}, \\
     C = \begin{bmatrix*}[r] % matrix C
       -0.92 & 0.94 \\
       -0.08 & 0.33 \\
       0.43 & -0.70 \\
       0.19 & -0.74
       \end{bmatrix*}, \ & D = \begin{bmatrix*}[r] % matrix D
         0.82 & -0.49 & 0.53 & 0.27 \\
         0.84 & 0.88 & -0.68 & -0.94 \\
         -0.72 & -0.20 & 0.70 & 0.26 \\
         -0.90 & -0.23 & -0.33 & 0.56
         \end{bmatrix*}.    
   \end{aligned}
 \end{equation*} 
For this system,
 the dual LMI \eqref{eq:RedLMI_DD} turns out to be feasible, 
 and the resulting dual solution $H$ is numerically verified to 
 be $\rank{H} = 1$. 
The full-rank factorization of $H$ as well as
 $\z,\w\in\mb{R}^4$ in \textit{Theorem \ref{thm:DD}} are obtained as
 \begin{equation*}
 \scalebox{0.7}
 {$
   h_1 = \begin{bmatrix*}[r]
     1.7238 \\
     -0.1691
     \end{bmatrix*}, \ h_2 = \begin{bmatrix*}[r]
       -0.2516 \\
       -0.0279 \\
       0.2033 \\
       0.0279
       \end{bmatrix*}, \ \z = \begin{bmatrix*}[r]
         -1.8222 \\
         -0.5940 \\
         1.1959 \\
         0.6340
         \end{bmatrix*}, \ \w = \begin{bmatrix*}[r]
           -0.2516 \\
           -0.0279 \\
           0.2033 \\
           0.0279
           \end{bmatrix*} (= h_2).
 $}          
 \end{equation*}
It is obvious that $h_1 \ne 0$ (the assertion (i) of \textit{Theorem \ref{thm:DD}}). 
In Fig.~\ref{fig:in-out_DD},
 the solid line represents the input-output characteristics of 
 the destabilizing nonlinear (piecewise linear) operator $\phiwc$ constructed by 
 following the procedure in the preceding subsection, 
 while the dashed line represents a line of slope one. 
From Fig.~\ref{fig:in-out_DD},
 we can readily see that $\phiwc$ is odd and $\slope{\phiwc}{0}{1}$ (the assertion (ii)). 
Fig.~\ref{fig:trajectory_DD} shows 
 the vector field of the system $\Sigma$ with the nonlinearity $\Phi = \Phiwc$, 
 together with the state trajectories from 
 initial states $x(0) = h_1$ and $x(0) = [-0.5\ -0.5]^T$. 
The state trajectory from 
 the initial state $x(0) = [-0.5\ -0.5]^T$ converges to the origin. 
However, as proved in \textit{Theorem \ref{thm:DD}},
 the state trajectory from the initial state $x(0) = h_1$ does not evolve and 
 $x(0) = h_1$ is confirmed to be an equilibrium point of 
 the feedback system $\Sigma$ (the assertion (iii)).
%%%%%%%%%%%%%%%%%%%%%%%%%%%%%%%%%%%%%%%%%%%%%%%%%%%%%%%%%%%%%%%%%%%%%%%%%%
%%%%%%%%%%%%%%%%%%%%%%%%%%%%%%%%%%%%%%%%%%%%%%%%%%%%%%%%%%%%%%%%%%%%%%%%%%
\section{Conclusion}
This paper presented novel results on the absolute stability analysis of
discrete-time nonlinear feedback systems.  
Namely,  if the solution of the proposed dual LMI satisfies a specific rank condition, 
then we showed that we can detect destabilizing nonlinear operator
that enables us to conclude that the feedback system of interest is 
never absolutely stable.  
We illustrated the technical results by numerical examples.  

In the present paper, the main results have been obtained for 
slope-restricted nonlinearities of $\mathrm{slope}[0,1]$.  
Therefore, it is an important future topic to 
extend the results to the general case $\mathrm{slope}[\mu,\nu]\ (\mu \le 0 \le \nu)$.  
In addition, we note that the recent study (\cite{Yuno_CDC2024}) has proposed
a new type of static OZF multipliers for slope-restricted and {\it idempotent} nonlinearities. 
Therefore, expanding the scope of the current results to
the set of new multipliers is also a challenging future topic.
These topics are currently under investigation.

% \begin{ack}
% Place acknowledgments here.
% \end{ack}

%\bibliography{rocond2025}             % bib file to produce the bibliography

                                                     % with bibtex (preferred)
            
% \appendix
% \section{A summary of Latin grammar}    % Each appendix must have a short title.
% \section{Some Latin vocabulary}              % Sections and subsections are supported  
                                                                         % in the appendices.
\end{document}